\documentclass[conference]{IEEEtran}
\IEEEoverridecommandlockouts
\usepackage{cite}
\usepackage{amsmath,amssymb,amsfonts}
\usepackage{algorithmic}
\usepackage{graphicx}
\usepackage{epsfig} 
\usepackage{textcomp}
\usepackage{bm}
\usepackage{xcolor}
\usepackage{times}
\usepackage{tabularx}
\usepackage{multirow}
\usepackage{booktabs} 
\usepackage{subcaption}
\usepackage{enumitem}
\usepackage[labelsep=period, font=small,labelfont=bf]{caption}

\def\BibTeX{{\rm B\kern-.05em{\sc i\kern-.025em b}\kern-.08em
    T\kern-.1667em\lower.7ex\hbox{E}\kern-.125emX}}

\begin{document}

\title{Energy Storage Arbitrage Under Price Uncertainty: Market Risks and Opportunities
}

\author{\IEEEauthorblockN{ Yiqian Wu}
\IEEEauthorblockA{\textit{Dept. of Electrical Engineering} \\
\textit{Columbia University}\\
New York, NY 10027, USA \\
\texttt{yiqian.wu2@columbia.edu}}
\and
\IEEEauthorblockN{Bolun Xu}
\IEEEauthorblockA{\textit{Dept. of Earth and Environmental Engineering} \\
\textit{Columbia University}\\
New York, NY 10027, USA \\
\texttt{bx2177@columbia.edu}}
\and
\IEEEauthorblockN{James Anderson}
\IEEEauthorblockA{\textit{Dept. of Electrical Engineering} \\
\textit{Columbia University}\\
New York, NY 10027, USA \\
\texttt{james.anderson@columbia.edu}}
}

\maketitle

\begin{abstract}
 We investigate the profitability and risk of energy storage arbitrage in electricity markets under price uncertainty, exploring both robust and chance-constrained optimization approaches. We analyze various uncertainty representations, including polyhedral, ellipsoidal uncertainty sets and probabilistic approximations, to model price fluctuations and construct efficient frontiers that highlight the tradeoff between risk and profit. Using historical electricity price data, we quantify the impact of uncertainty on arbitrage strategies and compare their performance under distinct market conditions. 
  The results reveal that arbitrage strategies under uncertainties can effectively secure expected profits, and robust strategies perform better in risk management across varying levels of conservativeness, especially under highly volatile market conditions.
  This work provides insights into storage arbitrage strategy selection for market participants with differing risk preferences, emphasizing the adaptability of efficient frontiers to the electricity market.
\end{abstract}

\section{Introduction}
Increasing intermittent renewable resources presents significant challenges to grid operation, and energy storage systems are essential for balancing supply and demand. Energy storage participants in electricity markets leverage price volatility to arbitrage price differences based on forecasts of future prices, making a profit while aiding grid operations to reduce peak demands. 
However, with the increasing complexity of the power grid, the uncertainty in price forecasting has also inevitably grown. 
Given that storage profits are highly sensitive to price fluctuations, it is essential to understand the implications of uncertainty on profitability to develop resilient and effective strategies for market participation.

Stochastic optimization is commonly used to manage storage arbitrage decisions under price prediction uncertainties~\cite{zheng2022arbitraging, krishnamurthy2017energy}. However, it requires substantial computational resources to address uncertainty representations, which grow exponentially with longer time horizons and often lack a performance guarantee. For price arbitrage in particular, storage operators face the risk of negative returns if they charge high prices and are unable to sell the energy at profitable times. While regular market participants can absorb occasional losses, this risk significantly deters participation from non-professional players, like behind-the-meter home batteries, vehicle-to-grid, or community storage systems. This is especially true during grid contingencies when price volatility is high, and the potential for loss is amplified.


Risk management strategies for storage market participants have been explored through Conditional Value at Risk (CVaR)\cite{dicorato2009}, where Monte Carlo simulations are used to illustrate efficient frontiers that capture the tradeoff between profit and CVaR. Multi-level optimization approaches with robust objectives have also been studied to address uncertainties in rivals' offers within auction-based markets, presenting both worst-case strategic net revenue and expected profit-based efficient frontiers across varying levels of conservativeness\cite{fanzeres2019, mortaz2020}. More extensive reviews of related work are presented in~\cite{roald2023}. However, these methods still require substantial computational power to generate scenarios and simulate grid operations.

This paper proposes a computationally-efficient risk-averse arbitrage framework for energy storage. This framework is especially suitable for non-professional storage to arbitrage with controlled risk based on the unit's availability occasionally. Our paper makes the following contributions:
\begin{enumerate}
    \item We propose a general uncertainty-incorporated storage arbitrage formulation that can accommodate a variety of price uncertainty models and risk preferences.
    \item We present a performance evaluation framework via efficient frontiers highlighting the tradeoff between risk and profit.
    \item We compare mainstream strategies under uncertainty and provide insights into risk strategy selection.
\end{enumerate}
The paper is organized as follows: Section II introduces the formulation, and Section III presents the uncertainty representation. Section IV discusses the strategy selection, Section V presents the case study, and Section IV concludes the paper.

\section{Problem Formulation}

\subsection{Energy Storage Arbitrage}
We formulate energy storage arbitrage as a multi-interval self-scheduling problem according to a series of price forecasts $\lambda_t$ for all $t \in \mathcal T$, where $\mathcal T = \{1,2,\dots, T\}$. Denote the price series  as \(\lambda = \{\lambda_1, \lambda_2, \dots, \lambda_T\}\). 
Assuming the energy storage acts as a price taker in the electricity market, meaning its operation does not influence market clearing prices, the arbitrage problem with perfect price forecasting is modeled by the linear program~\cite{wu2024}:
\allowdisplaybreaks
\begin{subequations}
  \label{eq:sses_taker}
\begin{align}
  \underset{p_t,b_t, e_t}{\text{maximize}}&\quad \textstyle\sum_{t \in \mathcal{T}}\lambda_t (p_t-b_t)  \notag \\
  \mathrm{s.t.}
  &\quad 0\le p_t,b_t \le P \label{eq:pch_limit} \\
  & \quad p_t = 0 \textrm{ if } \lambda_t <0 \label{eq:storage_idle}\\
  &\quad e_t - e_{t-1} = -{p_t}/{\eta} + b_t \eta  \label{eq:soc_dynamic}\\
  & \quad 0\le e_t \le E \label{eq:soc_bar} 
\end{align}
\end{subequations}
where the operation decisions include the discharging, $p_t$, and charging, $b_t$, power over time period $t$. The objective function describes the expected arbitrage profit from storage decisions (we omit operational costs here for simplicity, thus profit equals to revenue).
Constraint~\eqref{eq:pch_limit} imposes the power rating limits for discharging and charging. 
Constraint~\eqref{eq:storage_idle} provides a convex relaxation to ensure the storage does not discharge during periods of negative pricing, which serves as a sufficient condition to prevent simultaneous discharging and charging~\cite{xu2020}.
The inter-temporal dynamics of the state of charge (SoC) $e_t$ are captured in~\eqref{eq:soc_dynamic}, taking into account discharging and charging efficiencies $\eta \in (0,1]$. Inequality~\eqref{eq:soc_bar} models the energy storage capacity.
We have normalized the time period duration into $p_t$ and $b_t$, so no duration coefficient is needed.
We denote by $\Omega$ the feasible set of~\eqref{eq:sses_taker}.

The price signal $\lambda$ in the objective function represents the forecasted electricity prices at the time of storage decision-making. As a result, arbitrage decisions are subject to forecasting uncertainty, potentially leading to suboptimal outcomes and risks. 
The remainder of this section explores different modeling approaches to characterize price uncertainty and further address decision-making under uncertainty.


\subsection{Robust Decision Making}
Robust optimization is a deterministic worst-case approach to handling uncertainty. Each point in the uncertainty set is equally likely to occur and the optimization problem is formulated to maximize the performance of the worst-case realization.
We assume that electricity price forecasts fall within nonempty bounded uncertainty sets, i.e., $\lambda_t \in \mathcal{U}_t $, where $\mathcal{U}_t$ represents the uncertainty set for electricity prices $\lambda_t$, and define $\mathcal{U}:=\bigcup \mathcal{U}_t$.
The robust counterpart of the arbitrage model~\eqref{eq:sses_taker} is then as:
\begin{align}
  \label{eq:taker_ro}
  \underset{\{p_t,b_t,e_t\} \in \Omega, \gamma}{\text{maximize}} \  \gamma \quad
  \mathrm{s.t.} \ \underset{\lambda_t \in \mathcal{U}_t}{\text{minimize}}  \ \textstyle\sum_{t \in \mathcal{T}}\lambda_t (p_t-b_t) \geq \gamma.
\end{align}

One of the main challenges in the robust optimization approach is the modeling of uncertainty sets, which directly impacts both the efficiency and robustness of decision-making. Commonly used uncertainty sets include the box set, polyhedral sets~\cite{korolko2017}, and ellipsoidal sets~\cite{golestaneh2018}. 
Box sets are the simplest to construct, but as they lack the ability to account for correlation, are often too conservative in practice. Polyhedral sets offer the most feasibility but are difficult to fit to data. Ellipsoidal sets provide a tractable tradeoff between the two.

\subsection{Chance-Constrained Decision Making}

In certain special cases, it is feasible to work directly with probabilistic uncertainty models when the uncertainty can be well-approximated by specific tractable distributions. Assuming these conditions are met, 
chance constraints enable a risk-averse solution within the stochastic framework by limiting the probability of constraint violations:
    \begin{align}
   \underset{\{p_t,b_t,e_t\} \in \Omega, \gamma}{\text{maximize}} \ \gamma \qquad
      \mathrm{s.t.} \ \mathbb{P}( \textstyle\sum_{t \in \mathcal{T}}\lambda_t (p_t-b_t) \geq \gamma) \ge \Gamma .\label{eq:prob_con}
    \end{align}
where the chance constraint in problem~\eqref{eq:prob_con} ensures that the stochastic objective reaches  optimality with a probability (or confidence) exceeding $\Gamma$.

\section{Uncertainty Representation}
We consider three main uncertainty representations: 1) polyhedral, 2) ellipsoidal, and 3) probabilistic. For each category, we introduce two specific formulations, each corresponding to a different arbitrage strategy, to characterize and manage the uncertainty in electricity prices. 
Let $\Gamma$ denote the user-defined uncertainty budget, indicating the level of conservativeness in the strategy. A higher $\Gamma$ widens the coverage, i.e., captures more price uncertainty, leading to more conservative decision-making. Conversely, lower values of  $\Gamma$ result in more aggressive strategies.

\subsection{Polyhedral uncertainty set}
The polyhedral uncertainty set is defined as:
\begin{align}
  \mathcal{U} = \left\{ \lambda \in \mathbb{R}^T \mid D\lambda \leq d \right\} \notag
\end{align}
where $D \in \mathbb{R}^{m \times T}$ and $d \in \mathbb{R}^m$ are fixed.

Given the uncertainty set defined above, deploying strong duality theorem of linear programs, the constraints encoding $\mathcal U$ in problem~\eqref{eq:taker_ro} yield the following formulation~\cite{korolko2017}:
  \begin{align}
      -y^{\mathsf{T}} d   \ge \gamma,\quad  -y^{\mathsf{T}} D = {(p-b)}^{\mathsf{T}}  \notag
  \end{align}
where \(p = \{p_1, p_2, \dots, p_T\}\) and \(b = \{b_1, b_2, \dots, b_T\}\) denote the sequence of storage discharge and charge actions, and $y$ represents the vector of dual variables.

Some representative examples of the polyhedral uncertainty sets are listed below:
    \begin{enumerate}[leftmargin=*]           \item \textbf{Poly-Quantile:} 
        \( \mathcal{U}(\Gamma) := \{ \lambda \in \mathbb{R}^{T}  \mid  z_{1-\Gamma }(\lambda_t)  \leq \lambda_t \leq  z_{\Gamma }(\lambda_t) , \ \forall t \in \mathcal T \} \), where \(  z_{\Gamma / 2}(\lambda_t) \) is the quantile of level $\Gamma/2$ of the electricity price at time $t$ estimated from historical data.

          \item \textbf{Poly-Mean-Std:} $ \mathcal{U}(\Gamma) := \{ \lambda \in \mathbb{R}^{T} \mid \left| \frac{\lambda_t - \mathbb{E} \lambda_t}{\sigma(\lambda_t)} \right| \leq \Gamma, \ \forall t \in \mathcal T \} $,
           where \( \sigma(\lambda_t) \) is the standard deviation based on historical observations of electricity price  at time $t$. 

    \end{enumerate}

    \subsection{Ellipsoidal uncertainty set}

    The ellipsoidal uncertainty set is defined as:
    \begin{align}
    \mathcal{U} = 
    \left\{ \bar{\lambda} + Q u  \in \mathbb{R}^T  \mid \| u \| \leq 1 \right\} \notag 
    \end{align}
    where $Q\in \mathbb{S}_{+}^{T }$ represents a fixed $T\times T$  positive semidefinite matrix, and \( \bar{\lambda } \) is the nominal point, i.e., the ellipse center.
    Note that if \( Q \) is not full-rank, it means there is no uncertainty in a given coefficient of~\( \lambda \).

   Accordingly, problem~\eqref{eq:taker_ro} can be rewritten as a second-order cone program (SOCP)~\cite{boyd2004} and can be efficiently solved using off-the-shelf convex optimization solvers. The robust constraint in problem~\eqref{eq:taker_ro} is recast as:
      \begin{align}
        \bar{\lambda}^{\mathsf{T}} (p-b) -   \|Q^{\mathsf{T}} (p-b) \|   \ge \gamma .\notag
      \end{align}
    
    Some representative examples of the ellipsoidal uncertainty sets are listed below:
    \begin{enumerate}
      \item \textbf{Ellip-Min-Vol:} 
    $Q$ can be obtained by fitting the historical data of electricity prices into a minimum volume enclosing ellipsoid via the following model:
    \begin{subequations}
\begin{align*}
    \textstyle \min_{Q\in \mathbb{S}_{++}^{T }} &\quad -\mathrm{\log\, \det}  \left(Q ^{-2}\right)\\
     \text{s.t.} &\quad {(\lambda^s - \bar{\lambda})}^{\mathsf{T}} Q^{-2} ( {\lambda}^s -\bar{\lambda}) \leq \Gamma^{-2},\, \forall s \in \mathcal{S} 
\end{align*}
\end{subequations}
where $\lambda^s$ is the $s$-th sample of the electricity prices, $\mathcal S$ is the set of historical prices.

\item \textbf{Ellip-Cov:} 
\(
\mathcal U (\Gamma) =\{ \lambda \in \mathbb{R}^{T} \mid {(\lambda - \bar{\lambda})}^{\mathsf{T}} \Sigma^{-1} (\lambda - \bar{\lambda}) \leq \Gamma^{-2} \} 
\),
where $\Sigma$ is the covariance matrix of electricity prices estimated from historical data, \(
  \Sigma_{ij} = \frac{1}{m - 1} \sum_{k=1}^{m} (x_{ki} - \bar{x}_i)(x_{kj} - \bar{x}_j)
  \), and $Q^{\mathsf{T}}Q = \Sigma$.


\end{enumerate}

\subsection{Probability distribution approximation}
The assumption of a normal distribution for electricity prices is widely adopted in the literature~\cite{xu2020}, while others suggest a lognormal distribution~\cite{conejo2002}. In our case studies, we consider both distributions for the probability approximation: 
\begin{enumerate}
    \item \textbf{Chance-Normal:} assuming the price forecast following normal distributions, with mean $\mu_t$ and variance $\sigma_t^2$, the chance constraint in problem~\eqref{eq:prob_con} can be expressed as a second-order cone constraint~\cite{boyd2004}:
    \begin{align}
\textstyle \sum_{t \in \mathcal{T}} \mu_t (p_t - b_t) - \Phi^{-1}(\Gamma) \sqrt{\sum_{t \in \mathcal{T}} \sigma_t^2 {(p_t - b_t)}^2} \geq \gamma  \nonumber 
\end{align}
where $\Phi^{-1}(\cdot)$ is the inverse of the standard normal cumulative distribution function, $\Phi(z) = \frac{1}{\sqrt{2\pi}} \int_{-\infty}^{z} e^{-\tau^2 / 2} \, d\tau  \nonumber$. 

\item \textbf{Chance-LogNormal:} models the price forecast as lognormal distributions, i.e., $\lambda_t \sim LN(\mu_t, \sigma_t^2)$.
Since the left-hand side of the probability expression in the chance constraint involves a sum of lognormals, which is complex to handle directly, we can approximate it by a single lognormal variable, $ Y \sim \mathrm{LN}(\mu', \sigma'^2) $, using moment matching or more specialized methods, such as the Fenton-Wilkinson approximation~\cite{mehta2007}. This allows us to apply the inverse cumulative distribution function of the lognormal distribution, $F(y) = \Phi(\frac{\log y-\mu'}{\sigma'})$, similar to the approach for handling Chance-Normal.


\end{enumerate}

\section{Risk-Profit Analysis and Strategy Selection}
We use efficient frontier to exhibit the correlation between risk and profit to evaluate the quality of the strategies and decisions in the context of energy storage arbitrage. 
More specifically, we define the key metrics for analysis as follows:

\noindent\textbf{Expected profit:} The profit or loss that a market participant anticipates based on historical records. While expected profit is not guaranteed, historical data provides a reasonable basis for forecasting. Therefore, the expected profit can be considered a long-term weighted average of historical profits. In our study, it is equivalent to the averaged daily profit in the out-of-sample experiments.

\noindent\textbf{Risk:} The number of days with negative profits per year in the out-of-sample experiments.

The tradeoff between risk and profit is balanced via the uncertainty budget $\Gamma$. The storage owner solves the arbitrage problem for different values of the parameter $\Gamma$, and creates the efficient frontier in the profit-risk space. 
This efficient frontier then serves as a decision-making framework, allowing the storage owner to select strategies aligned with their risk preferences and profit objectives.
Additionally, we calculate the ratio of non-negative profit days to the total number of samples as another performance metric.


\begin{figure}[!t]
    \centering
    \includegraphics[width=0.45\textwidth]{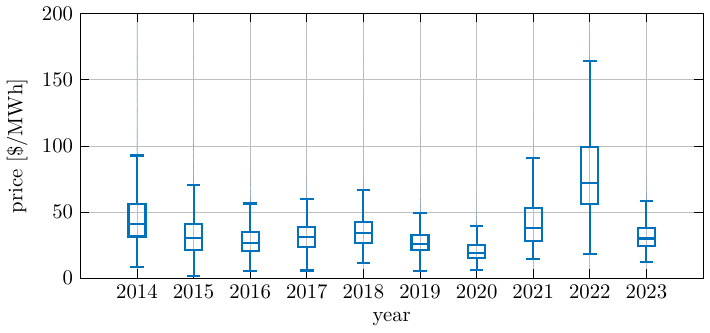}
    \vspace{-0.6em}
    \caption{NYC LMP from 2014 to 2023 (extreme outliers omitted).}
    \label{fig:price}
    \vspace{-0.5em}
\end{figure}

\section{CASE STUDIES}
We examine the profit performance of different energy storage arbitrage approaches under uncertainty by simulating a 24-interval dispatch in the day-ahead market. Our simulation features a 2.5MW/10MWh storage unit that starts and ends at 50\% SoC, with discharging and charging efficiencies of 0.9. 
We perform case studies through out-of-sample arbitrage, focusing on three key elements: 1) market conditions, evaluated by applying arbitrage decisions to different pricing years, 2) uncertainty representation, including the choice of uncertainty set forms and probability distribution approximations, and 3) dataset selection, used to quantify price uncertainty. In all figures, overlapping curves are slightly offset for clarity.

\subsection{Data Description}
We use hourly price data collected from the New York ISO (NYISO) for the New York City (NYC) zone from 2014 to 2023 to estimate uncertainty characteristics~\cite{NYISO_price}. The price dataset statistics are visualized in Fig.~\ref{fig:price}. 
Prices generally range between \$0/MWh and \$200/MWh, displaying a relatively consistent pattern over the ten-year period. 
Annual variations can be attributed to fluctuating natural gas prices, transmission congestion conditions, global market dynamics, etc.~\cite{patton2024}, allowing us to evaluate the robustness and profitability of different arbitrage strategies under various conditions.

Uncertainty quantification is conducted using price data from 2014 to 2021, with the last two years (2022-2023) reserved for out-of-sample testing, unless otherwise specified. The sample characteristics are as follows: 1) in 2022, prices reached higher levels with significantly increased volatility, resulting in a pronounced ``fat tail'' in the distribution for that year, and 2) in 2023, prices remained at more stable levels and aligned more closely with the 2014-2021 averages.

\subsection{Arbitrage under Uncertainty}
The 24-hour ahead arbitrage problem is solved for different values of the uncertainty budget~$\Gamma$, with 2022 as the test year, in which prices exhibit high volatility. 
We use the historical mean of electricity prices from 2014 to 2021 as the nominal ellipse center $\bar{\lambda}$. The worst-case profit, expected profit, and non-negative profit ratio are calculated for each strategy. The results are summarized in Table~\ref{table:summary}\footnote{$\Gamma$ is normalized and constrained within the range $\Gamma \in [0,1]$ to simplify the presentation of simulation results. Here, $\Gamma = 0$ corresponds to the most aggressive strategy, while $\Gamma = 1$ represents the highest level of conservativeness, beyond which no profit is generated (since chance-constrained strategies never yield zero profit, $\Gamma = 1$ represents the most conservative scenario).}.

\begin{table}[t]
    \centering
    \caption{Arbitrage Performance under Uncertainty}
    \label{table:summary}
\vspace{-0.5em}
\begin{subtable}[t]{0.49\textwidth}
    \centering
    \scriptsize
    \caption{Worst-case profit}
    \vspace{-0.5em}
    \label{table:summary_1}
    \begin{tabular}{
        p{0.09\columnwidth}|
        p{0.08\columnwidth}p{0.08\columnwidth}p{0.08\columnwidth}p{0.08\columnwidth}p{0.08\columnwidth}p{0.13\columnwidth}}
        \toprule
                    \begin{tabular}[c]{@{}l@{}}\textbf{Strategy}\\\end{tabular} & \begin{tabular}[c]{@{}l@{}}\textbf{Poly}\\ \textbf{Quantile}\end{tabular} & \begin{tabular}[c]{@{}l@{}}\textbf{Poly}\\ \textbf{Mean-Std}\end{tabular} & \begin{tabular}[c]{@{}l@{}}\textbf{Ellip}\\ \textbf{Min-Vol}\end{tabular} & \begin{tabular}[c]{@{}l@{}}\textbf{Ellip}\\ \textbf{Cov}\end{tabular} & \begin{tabular}[c]{@{}l@{}}\textbf{Chance}\\ \textbf{Normal}\end{tabular} & \begin{tabular}[c]{@{}l@{}}\textbf{Chance}\\ \textbf{LogNormal}\end{tabular} \\
                   \midrule
$\Gamma$ = 0   & 84.00               & 103.13             & 103.13  & 103.13 & 101.15 & 99.03     \\
$\Gamma$ = 0.2 & 61.06               & 71.73              & 72.14   & 73.61  & 82.11  & 87.37     \\
$\Gamma$ = 0.4 & 42.51               & 47.14              & 42.41   & 47.68  & 62.99  & 74.88     \\
$\Gamma$ = 0.6 & 27.86               & 28.81              & 18.75   & 24.62  & 52.44  & 61.30     \\
$\Gamma$ = 0.8 & 13.08               & 11.78              & 7.03    & 8.82   & 52.44  & 47.20     \\
$\Gamma$ = 1   & 0.00                & 0.00               & 0.00    & 0.00   & 52.44  & 45.73     \\
    \bottomrule
\end{tabular}
    \vspace{0.3em}
\end{subtable}
\hspace{0.05\textwidth} 
    \vspace{-0.3em}
\begin{subtable}[t]{0.49\textwidth}
    \centering
    \caption{Expected profit}
    \label{table:summary_2}
        \vspace{-0.3em}
    \scriptsize
    \begin{tabular}{
        p{0.09\columnwidth}|
        p{0.08\columnwidth}p{0.08\columnwidth}p{0.08\columnwidth}p{0.08\columnwidth}p{0.08\columnwidth}p{0.13\columnwidth}}
        \toprule
    \begin{tabular}[c]{@{}l@{}}\textbf{Strategy}\\\end{tabular} & \begin{tabular}[c]{@{}l@{}}\textbf{Poly}\\ \textbf{Quantile}\end{tabular} & \begin{tabular}[c]{@{}l@{}}\textbf{Poly}\\ \textbf{Mean-Std}\end{tabular} & \begin{tabular}[c]{@{}l@{}}\textbf{Ellip}\\ \textbf{Min-Vol}\end{tabular} & \begin{tabular}[c]{@{}l@{}}\textbf{Ellip}\\ \textbf{Cov}\end{tabular} & \begin{tabular}[c]{@{}l@{}}\textbf{Chance}\\ \textbf{Normal}\end{tabular} & \begin{tabular}[c]{@{}l@{}}\textbf{Chance}\\ \textbf{LogNormal}\end{tabular} \\
                   \midrule
$\Gamma$ = 0   & 194.99              & 194.99             & 194.99  & 194.99 & 194.99 & 194.99    \\
$\Gamma$ = 0.2 & 182.93              & 182.93             & 192.99  & 191.93 & 192.99 & 194.99    \\
$\Gamma$ = 0.4 & 182.93              & 164.10             & 190.73  & 187.12 & 175.70 & 194.99    \\
$\Gamma$ = 0.6 & 154.22              & 154.22             & 116.73  & 157.26 & 67.56  & 165.48    \\
$\Gamma$ = 0.8 & 154.22              & 154.22             & 59.56   & 131.28 & 67.56  & 131.64    \\
$\Gamma$ = 1   & 0.00                & 0.00               & 0.00    & 0.00   & 67.56  & 53.74     \\
    \bottomrule
\end{tabular}
        \vspace{-0.3em}
\end{subtable}
\hspace{0.05\textwidth} 

\begin{subtable}[t]{0.49\textwidth}
    \centering
    \caption{Non-negative profit ratio}
            \vspace{-0.3em}
    \label{table:summary_3}
    \scriptsize
    \begin{tabular}{
        p{0.09\columnwidth}|
        p{0.08\columnwidth}p{0.08\columnwidth}p{0.08\columnwidth}p{0.08\columnwidth}p{0.08\columnwidth}p{0.13\columnwidth}}
        \toprule
                \begin{tabular}[c]{@{}l@{}}\textbf{Strategy}\\\end{tabular}  & \begin{tabular}[c]{@{}l@{}}\textbf{Poly}\\ \textbf{Quantile}\end{tabular} & \begin{tabular}[c]{@{}l@{}}\textbf{Poly}\\ \textbf{Mean-Std}\end{tabular} & \begin{tabular}[c]{@{}l@{}}\textbf{Ellip}\\ \textbf{Min-Vol}\end{tabular} & \begin{tabular}[c]{@{}l@{}}\textbf{Ellip}\\ \textbf{Cov}\end{tabular} & \begin{tabular}[c]{@{}l@{}}\textbf{Chance}\\ \textbf{Normal}\end{tabular} & \begin{tabular}[c]{@{}l@{}}\textbf{Chance}\\ \textbf{LogNormal}\end{tabular} \\
                   \midrule
$\Gamma$ = 0   & 0.93                & 0.93               & 0.93    & 0.93   & 0.93   & 0.93      \\
$\Gamma$ = 0.2 & 0.96                & 0.96               & 0.90    & 0.94   & 0.90   & 0.93      \\
$\Gamma$ = 0.4 & 0.96                & 0.94               & 0.89    & 0.96   & 0.85   & 0.93      \\
$\Gamma$ = 0.6 & 0.94                & 0.94               & 0.82    & 0.95   & 0.57   & 0.91      \\
$\Gamma$ = 0.8 & 0.94                & 0.94               & 0.70    & 0.93   & 0.57   & 0.89      \\
$\Gamma$ = 1   & 1.00                & 1.00               & 1.00    & 1.00   & 0.57   & 0.69                               \\
    \bottomrule                      
    \end{tabular}
 
    \end{subtable}
               \vspace{-1.5em}
\end{table}

As conservativeness increases, both the worst-case profit and out-of-sample expected profit decrease, as anticipated -- this tradeoff can be referred to as \emph{the price of robustness}. 
Different from the robust strategies (the first four strategies), chance-constrained strategies guarantee positive profits even in the most conservative scenarios, which is attributable to the probabilistic nature of the constraints, and level off as $\Gamma$ exceeds 0.8.
As for the expected profit, within a moderate 0-0.4 range of $\Gamma$, the expected profit remains relatively stable, with a slight decrease (within 15\%) as $\Gamma$ increases. However, when $\Gamma$ reaches 0.6, the expected profit drops significantly, indicating that the strategy becomes overly conservative, leading to a loss of profit opportunities.
When $\Gamma$ is set to 1.0, the expected profit is reduced to zero, indicating that the strategy is too conservative to yield any profit.
For chance-constrained strategies, however, the expected profit remains positive, similar to the behavior of the worst-case profit.


\begin{figure}[!t]
  \centering
  \begin{subfigure}[b]{0.45\textwidth}
      \centering
      \includegraphics[width=\textwidth]{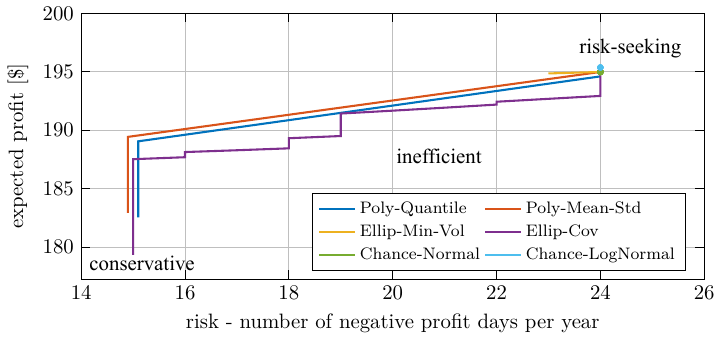}
                 \vspace{-1.7em}
      \caption{Tested with prices from year 2022}
  \end{subfigure}
  
  \begin{subfigure}[b]{0.45\textwidth}
      \centering
      \includegraphics[width=\textwidth]{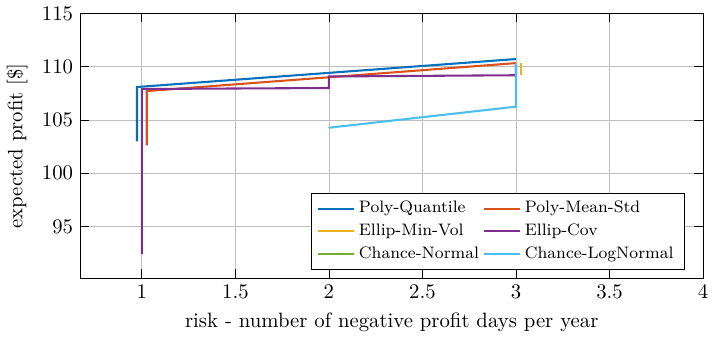}
                 \vspace{-1.7em}
      \caption{Tested with prices from year 2023}
  \end{subfigure}
  \caption{Efficient frontiers for different arbitrage strategies under uncertainty given variant market conditions.}
  \label{fig:diff_strategy}
\end{figure}

\begin{figure}[!t]
\centering
\begin{subfigure}[b]{0.45\textwidth}
    \centering
    \includegraphics[width=\textwidth]{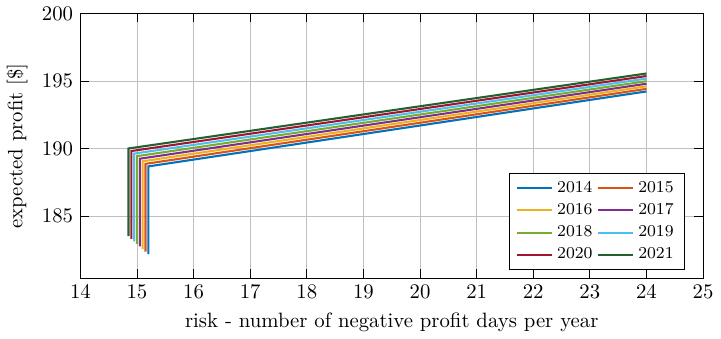}
                     \vspace{-1.7em}
    \caption{Poly-Mean-Std}
\end{subfigure}

\begin{subfigure}[b]{0.45\textwidth}
    \centering
    \includegraphics[width=\textwidth]{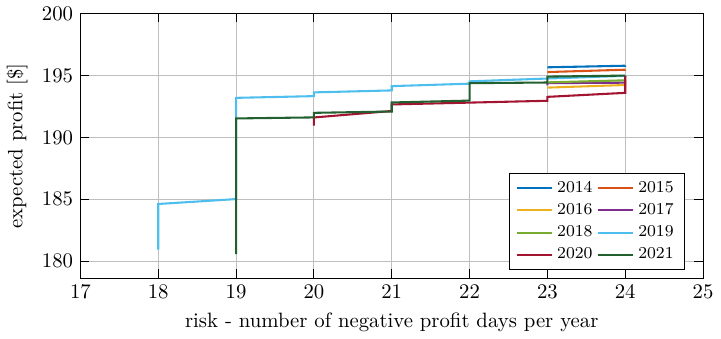}
                     \vspace{-1.7em}
    \caption{Ellip-Min-Vol}
\end{subfigure}

\begin{subfigure}[b]{0.45\textwidth}
    \centering
    \includegraphics[width=\textwidth]{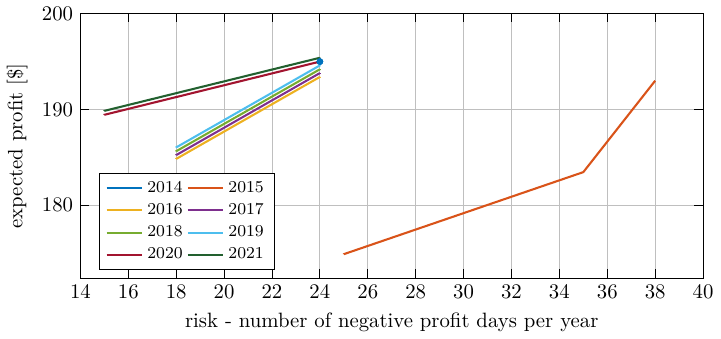}
                     \vspace{-1.7em}
    \caption{Chance-LogNormal}
\end{subfigure}
\caption{Efficient frontiers tested with prices from year 2022 for strategies characterized by datasets from year $n$ to 2021.}
\label{fig:diff_starting_year}                   
\end{figure}

The non-negative profit ratio remains close to 1.00 even under the most conservative scenarios, especially for polyhedral robust strategies. 
Ellip-Min-Vol and Chance-Normal strategies exhibit relatively weaker risk control, as the ratio dropping to 0.70 and 0.57, respectively,  while risk remains moderate for the other strategies as $\Gamma$  increases to 0.8. This may be due to the following reasons: 1) Ellip-Min-Vol aims to minimize the volume of the ellipsoid to reduce conservativeness, potentially compromising its effectiveness in capturing data correlations, and 2) Chance-Normal relies on a normal distribution approximation, which may not accurately represent the underlying price distribution.
One notable trend is that at lower levels of conservativeness, as $\Gamma$ increase, the non-negative profit ratio also increases, indicating that the strategies become more conservative and less risky. 
This trend captures a low-profit low-risk / high-profit high-risk shape of the efficient frontier -- a point we will discuss further in the following sections.


\subsection{Risk-Profit Analysis given Different Market Conditions}
Fig.~\ref{fig:diff_strategy} depicts efficient frontiers for arbitrage strategies under uncertainty tested on electricity price datasets from different years.
The frontier of year 2022 lies in the high-profit high-risk corner, indicating that the strategies are more profitable but riskier due to the higher average and greater volatility of electricity prices. 
In contrast, the frontier of year 2023 is at the low-profit low-risk end, reflecting reduced profitability as a result of lower average and volatility of electricity prices. These efficient frontiers effectively capture market conditions, revealing potential risks and opportunities and enabling market participants to select strategies that align with their risk preferences and profit objectives.

\subsection{Risk-Profit Analysis given Different Representations}
Fig.~\ref{fig:diff_strategy} also presents the performance of different strategies. The polyhedral strategies and Ellip-Cov covers a broader range for the efficient frontiers for both the 2022 and 2023 datasets, particularly in volatile markets, compared to other strategies.
This broader range allows potential risk reduction without compromising profit, demonstrating better performance in risk management. 
This is attributed to the fact that: 1) these uncertainty representations can effectively capture the strong correlations among uncertain parameters, and 2) chance-constrained strategies suffer from the drawback that they limit the probability of constraint violation yet fail to bound the magnitude, leading to a less effective risk management.

\subsection{Risk-Profit Analysis given Different Price Dataset}
Fig.~\ref{fig:diff_starting_year} depicts efficient frontiers for different strategies provided different historical datasets tested with prices from year~2022.
It shows that the ellipsoidal strategies exhibit the highest sensitivity to the historical data used, while polyhedral and chance-constrained strategies demonstrate greater robustness across datasets. This suggests that polyhedral strategies are particularly resilient to fluctuations in historical data.


\section{CONCLUSIONS}
We find that robust strategies outperform chance-constrained ones in risk management, particularly in volatile markets, while offering similar profit expectations. Notably, robust strategies using these representations capture correlations in price data, expanding the range of the efficient frontier and offering both high- and low-risk profit opportunities. The findings underscore the importance of selecting strategies that align with market participants' risk preferences. The framework provides a structured decision-making tool that adapts to evolving market conditions, offering practical implications for energy storage owners, especially for deploying risk-averse solutions in distributed storage solutions.
\bibliographystyle{IEEEtran}
\bibliography{uncertainty.bib}

\end{document}